\documentstyle[12pt]{article}
\titlepage
\newcommand{\ignore}[1]{}{}
\newcommand{\f}{\frac}

\renewcommand{\(}{\left(}
\renewcommand{\)}{\right)}
\renewcommand{\[}{\left[}
\renewcommand{\]}{\right]}

\newcommand{\no}{\nonumber}

\newcommand{\ep}{\epsilon}
\newcommand{\la}{\label}
\newcommand{\be}{\begin{eqnarray}}
\newcommand{\ee}{\end{eqnarray}}
\newcommand{\bestar}{\begin{eqnarray*}}
\newcommand{\eestar}{\end{eqnarray*}}

\newcommand{\si}{\sigma}

\newtheorem{th}{Theorem}
\newtheorem{lem}{Lemma}
\newtheorem{cor}{Corollary}

\begin{document}

\begin{center}

\vglue 1 cm

{\Large\bf Asymptotics of Studentized U-type processes
for changepoint problems}\\
[0.8cm]
{\large Mikl\'os Cs\"org\H{o}} \\
Carleton University, Ottawa, Canada
\\[0.5cm]
{\large Barbara Szyszkowicz} \\
Carleton University, Ottawa, Canada\\
[0.5cm]
{\large Qiying Wang} \\
University of Sydney, Australia
\end{center}

\bigskip
\centerline{\it Dedicated to the memory of Tibor Nemetz}

\bigskip
\centerline{{\bf ABSTRACT}}
\begin{center}
\begin{minipage}{120mm}

{This paper investigates weighted approximations for
 studentized $U$-statistics type processes,
 both with symmetric and antisymmetric kernels, only under the assumption
that the distribution of
the projection variate is in the domain of attraction of the normal law.
The classical second moment condition $E|h(X_1,X_2)|^2 < \infty$ is also relaxed in both cases.
The results can be used for  testing the null assumption of having a
random sample versus the alternative that there
 is a change in distribution in the sequence.
}

\end{minipage}
\end{center}



\medskip\noindent
{\bf Running Head:} Studentized U-type processes

\vspace {1cm}

\noindent
--------------------------------------------

\noindent
The research of M.~Cs\"org\H{o} and B.~Szyszkowicz is supported by their 
 NSERC Canada Discovery Grants at Carleton University, Ottawa, and Q.~Wang's research is
supported in part by Australian Research Council at University of Sydney.

\newpage


\section{Introduction and main results: the case of symmetric kernels
\la {sec1}}

Let $X, X_1, X_2,...$ be a sequence of non-degenerate independent real-valued random variables with distribution
function $F$. Suppose  we are interested in testing the   null hypothesis:
$$
H_0:\quad X_i, \, 1\le i\le n, \ have \ the \ same\ distribution,
$$
against the one change in distribution alternative:
\bestar
H_A: && there\ is\ an\ integer\ k,\ 1\le k<n, \ such\ that\\
&& P(X_1\le t)=\cdots= P(X_k\le t),\ P(X_{k+1}\le t)=\cdots=P(X_n\le t)\ \\
&&for \ all\ t\  and \
 P(X_k\le t_0)\neq P(X_{k+1}\le t_0)\  for\ some\ t_0.
\eestar

Testing for this kind of a change in distribution
  has been studied extensively  in the literature
by using  parametric as well as non-parametric methods.
One of the non-parametric methods was  proposed by Cs\"org\H{o} and
Horv\'ath (1988a, b),
who used  functionals of a $U$-statistics type
($U$-type, from now on) process to test $H_0$ against $H_A$.
Let $h(x,y)$ be a measurable real valued symmetric function, i.e. $h(x,y) = h(y,x)$.
The $U$-type process of Cs\"org\H{o} and Horv\'ath (1988a, b)
is defined by
$$
U_n(t)=Z_{[(n+1)t]}-n^2t(1-t)\theta, \qquad 0\le t\le 1,
$$
where $\theta=Eh(X_1, X_2)$, and
$$
Z_k=\sum_{i=1}^k\sum_{j=k+1}^nh(X_i,X_j),\qquad 1\le k < n.
$$

While $Z_k$ itself is not a $U$-statistic, it can be written as the sums of three $U$-statistics [cf. Cs\"org\H{o} and Horv\'ath (1988a, b, 1997)]. 
The rational behind the definition of $Z_k$ is comparing the first $k$ observations to the remaining $(n-k)$ ones for $k=1,\ldots,n-1$, via an appropriate bivariate kernel function $h(x,y)$ for the sake of capturing the possibility of having a change in distribution at an unknown time $k$ as postulated in $H_A$. Typical choices of symmetric kernel $h$ are $xy$, $(x-y)^2/2$ (the sample variance), $|x-y|$
(Gini's mean difference), and sign$(x+y)$ (Wilcoxon's one-sample statistic).

Throughout the paper,  we write
$
g(t)= E\(h(X, t)-\theta\),$  $ \sigma^2=Eg^2(X_1)$
and, for later use, we define a Gaussian process $\Gamma$ by
\be
\Gamma(t)=(1-t)\,W(t)\,+\,t\,\[W(1)-W(t)\], \quad 0\le t\le 1,
\la {ga1}
\ee
where $\{W(t), 0\le t<\infty\}$ is a standard Wiener process. Furthermore,
let $Q$ be the class of positive functions $q$ on $(0,1)$, i.e.,
$\inf_{\delta\le t\le 1-\delta}q(t)>0$ for $0<\delta<1$,
which are nondecreasing   in a neighbourhood of zero
and nonincreasing in a neighbourhood of one, and let
$$
I(q, c)=\int_{0+}^{1-}\f 1{t(1-t)}\exp\Big(-\f {cq^2(t)}{t(1-t)}\Big)dt,\qquad
0<c<\infty.
$$

In terms of these notations, Cs\"org\H{o} and Horv\'ath (1988a, b),
Szyszkowicz (1991, 1992) established  the following result
[cf.  Theorem 2.4.2 in Cs\"org\H{o} and Horv\'ath (1997)].

\medskip\noindent
\textbf{Theorem A}~ \textit{Assume $H_0$, $0<\sigma^2<\infty$
and $E|h(X_1, X_2)|^{2}<\infty$.
Then, on an appropriate probability space for $X,X_1, X_2,\cdots,$
we can define a sequence of Gaussian processes $\{\Gamma_n(t), 0\le t\le 1\}$
such that the equality in distribution
\be
\{\Gamma_n(t), 0\le t\le 1\}=_d
\{\Gamma(t), 0\le t\le 1\} \la {ad10}
\ee
holds for each $n\ge 1$, and as $n\to\infty$,
\be
\sup_{0 <t<1}\Big|\,n^{-3/2}\sigma^{-1}U_n(t)-
\Gamma_n(t)\,\Big|\Big/q(t)&=&o_P(1).
\la {th11}
\ee
if and only if $I(q,c)<\infty$ for all $c>0$.}

\medskip\noindent{\bf Remark 1}~
The condition $E|h(X_1,X_2)|^2 <\infty$ implies that $Eg^2(X_1)<\infty$, and we assume that $\sigma^2 = Eg^2(X_1)>0$. This is the so-called non-degenerate case when studying $U$-statistics via the function $g(t)=E(h(X,t)-\theta)$ that induces the projection of $U$-statistics into sums of i.i.d.\ random variables, the so-called Hoeffding (1948) projection principle that, in part, rests on a paper of Halmos (1946).


\medskip
For functions $x,\, y$ in $D[0,1]$ and $q\in Q$, we define the \textit{weighted sup-norm metric} $||/q||$ by
$$
||(x-y)/q|| = \sup_{0\leq t\leq 1} |(x(t)-y(t))/q(t)|,
$$
whenever this is well defined, i.e., {\it when} lim sup $|(x(t)-y(t))/q(t)|$ \textit{is finite} for $t\downarrow 0$ and $t\uparrow 1$.

In view of (2) and this terminology, (3) of Theorem A implies the following weak convergence, a functional limit theorem.

\medskip\noindent
\textbf{Corollary A} \textit{With $q\in Q$, and $\rightarrow_d$ standing for convergence in distribution as $n\to\infty$, we have}
$$
h\big\{n^{-3/2} \sigma^{-1} U_n(\cdot)/q(\cdot)\big\} \rightarrow_d h\{\Gamma(\cdot)/q(\cdot)\}
$$
\textit{for all} $h:D=D[0,1] \to I\!\!R$ {\it that are $(D,{\mathcal D})$ measurable and $||/q||$-continuous, or $||/q||$-continuous except at points forming a set of measure zero on $(D,{\mathcal D})$ with respect to the measure generated by the Gaussian  $\Gamma(\cdot)$ process, if and only if $I(q,c)<\infty$ for all $c>0$, where ${\mathcal D}$ denotes the $\sigma$-field of subsets of $D$ generated by the finite dimensional subsets of} $D$.

\medskip\noindent
\textbf{Remark A}~ For further use the statement of Corollary A will be summarized by writing, as $n\to\infty$,
$$
n^{-3/2} \sigma^{-1}U_n(\cdot)/q(\cdot) \Rightarrow \Gamma(\cdot)/q(\cdot) \hbox{ ~on~ } (D[0,1], {\mathcal D}, ||/q||).
$$

\medskip
For a summary of notions of convergence and weak convergence in general along these lines, we refer to pages 26--28 and Remarks 2 and 3 on page 49 of Shorack and Wellner (1986), and to Sections 3.3 and 3.4 of Cs\"org\H{o} (2002).

\medskip
Thus Theorem A provides a basic tool
for investigating the asymptotic behaviour of many test statistics for testing $H_0$ versus $H_A$ via corresponding
functionals of $\Gamma(\cdot)/q(\cdot)$ for appropriate choices of the kernel $h(x, y)$.
This, in turn,  motivates the establishment of our first result, in which
we  reduce the moment conditions related to the kernel $h(x, y)$.
It reads as follows.

\begin{th} \la {th1} Assume $H_0$, $0<\sigma^2<\infty$
and $E|h(X_1, X_2)|^{4/3}<\infty$.
Then, on an appropriate probability space for $X,X_1, X_2,\cdots,$
we can define a sequence of Gaussian processes $\{\Gamma_n(t), 0\le t\le 1\}$
such that {\rm(\ref {ad10})} holds true, and if $I(q,c)<\infty$ for some $c>0$,
then as $n\to\infty$,
\be
\sup_{1/n\le t\le (n-1)/n}\Big|\,n^{-3/2}\sigma^{-1}U_n(t)-
\Gamma_n(t)\,\Big|\Big/q(t)&=&o_P(1).
\la {th11aa}
\ee
\end{th}

In addition to  reducing the moment conditions required in Theorem A,
 the result (\ref {th11aa}) of Theorem \ref {th1} generalizes (\ref {th11}) as well.
 Namely,
 as a direct consequence of Theorem \ref {th1},
we have the following corollary.

\begin{cor} \la {cor1} Assume $H_0$, $0<\sigma^2<\infty$ and
$E|h(X_1, X_2)|^{4/3}<\infty$.
If $q\in Q$, then

\quad {\rm(a)}   we still have the conclusion of {\rm Theorem A}, i.e.,
{\rm(\ref {th11})} holds true
if and only if $I(q,c)<\infty$ for all $c>0$;

\quad {\rm(b)}  as $n\to\infty$,
\be
n^{-3/2}\sigma^{-1}U_n(\cdot)\,\Big/q(\cdot) &\Rightarrow&
\Gamma(\cdot)\,\Big/q(\cdot) \hbox{ ~on~ } (D[0,1], {\mathcal D}, ||/q||)
\la {adcor12}
\ee
if and only if $I(q,c)<\infty$ for all $c>0$;

\quad {\rm(c)} as $n\to\infty$,
\be
n^{-3/2}\sigma^{-1}\sup_{0<t<1}|U_n(t)|\Big/q(t)&\rightarrow_d&
\sup_{0<t<1}|\Gamma(t)|/q(t)
\la {th13}
\ee
if and only if $I(q,c)<\infty$ for some $c>0$.
\end{cor}

We note in passing that (a) implies (b), just like (3) implies Corollary A (cf.~(a) of Lemma 3).  However, (a) does not imply (c) (cf. (b) of Lemma 3).

\medskip

In view of the definition of $Z_k$, and hence also that of $U_n(t)$, when $\theta$ and $\sigma$ are known,
large values of the statistic on the left hand sides in
(\ref {th13}) for example, indicate a change in the distribution, and hence,
 based on Corollary~\ref {cor1}, rejection
of $H_0$ can be quantified accordingly.
Otherwise    $\theta$ and $\sigma$
need to be estimated.
A natural estimate of $\theta$ is
$$
\hat {\theta}=\frac 1{n(n-1)}\sum_{1\le i\ne j\le n}h(X_i, X_j),
$$
and that of $\sigma^2$ is
$$
\hat {\sigma}^2=\frac 1n \sum_{j=1}^n\Big (g(X_j)-\frac 1n\sum_{j=1}^ng(X_j)\Big)^2.
$$
According to the definition of $g(x)$,
$g(X_j)$ still depends on the usually unknown distribution function $F$ of $X$,
and hence it then can not be computed explicitly.
Since  we have that $g(x)+\theta=\int h(x,y)dF(y)$,
 we can replace $F$ by the empirical
distribution function $F_n$ of $X_1, X_2, ..., X_n$ under $H_0$. Consequently,
we may for example  estimate  $\sigma^2$ by
$$
\hat {\sigma}^2=\frac 1{n} \sum_{j=1}^n\Big (\f 1{n-1}\sum_{i=1\atop
i\ne j}^nh(X_i, X_j)
-\frac 1{n(n-1)}\sum_{1\le i\ne j\le n}h(X_i, X_j)\Big)^2.
$$
We note that this estimate  is in fact the jackknife estimator of $Var (\hat {\theta})$.
Now we may introduce a studentized U-type process as follows:
$$
\hat {U_n}(t)=n^{-3/2}(\hat {\sigma})^{-1}
\(Z_{[(n+1)t]}-n^2t(1-t)\hat {\theta}\), \qquad 0\le t\le 1.
$$
This process does not depend on the unknown parameters $\theta$ and $\si$ and
we now state the following main result of this paper, in which we replace the assumption that $0<\sigma^2<\infty$ by assuming only that $g(X_1)$ is in the domain of attraction of the normal law, written $g(X_1)\in \hbox{DAN}$ throughout.

\begin{th} \la {th2} Let $q\in Q$. Assume $H_0$, $E|h(X_1, X_2)|^{5/3}<\infty$
and that $g(X_1)\in \hbox{DAN}$. 
Then, on an appropriate probability space for $X,X_1, X_2,\cdots,$
we can define a sequence of Gaussian processes $\{\Gamma_n(t), 0\le t\le 1\}$
such that {\rm(\ref {ad10})} holds true and, as $n\to\infty$,
\be
\sup_{0< t<1}\Big|\,\hat {U_{n}}(t)-
\Gamma_n(t)\,\Big|\Big/q(t)&=&o_P(1),
\la {th15}
\ee
if and only if  $I(q,c)<\infty$ for all $c>0$. Consequently, as $n\to\infty$,
\be
\hat {U_{n}}(\cdot)/q(\cdot)&\Rightarrow&
\Gamma(\cdot)\,\Big/q(\cdot),\quad on\ (D[0,1], {\mathcal D}, ||/q||)
\la {adc12}
\ee
 if and only if
$I(q,c)<\infty$ for all $c>0$. Furthermore, as $n\to\infty$, we also have
\be
\sup_{0<t<1}|\hat {U_n}(t)|\Big/q(t)&\rightarrow_d& \sup_{0<t<1}|\Gamma(t)|/q(t)
\la {th17}
\ee
if and only if $I(q,c)<\infty$ for some $c>0$.
\end{th}

\noindent
\textbf{Remark 2}~ It is interesting and also of interest
to note that the class of  the weight functions in  (\ref {th17})
is bigger than that in (\ref {adc12}) [also compare  (\ref {th13})
with (\ref {adcor12})]. Such a phenomenon was first noticed
and proved for weighted empirical and quantile
processes by Cs\"org\H{o}, Cs\"org\H{o}, Horv\'ath and Mason [CsCsHM] (1986) and
then by Cs\"org\H{o} and Horv\'ath (1988b) for partial sums on assuming $E|X|^v<\infty$
for some $v>2$. For more details along these lines, we refer to
Szyszkowicz
 (1991, 1996, 1997),
and to Cs\"org\H{o}, Norvai\v{s}a and Szyszkowicz (1999).

\medskip\noindent
\textbf{Remark 3}~ As we noted already in Remark 1, the  condition that $0<\si^2=Eg^2(X_1)<\infty$ is the
so-called non-degenerate case when studying $U-$ statistics.
In Theorem \ref {th1} it  is a necessary  condition, while assuming
  $E|h(X_1, X_2)|^{4/3}<\infty$
is close to being necessary,   on account of the central
 limit theorem  for  $U$-statistics
(see Borovskikh (2002), for example).
Theorem \ref {th2}
puts a totally new countenance on the classical theory of
weak convergence for standardized
$U$-type process as in Theorem \ref {th1}
[cf. also Theorem A, Section 2.2.4 of  Cs\"org\H{o} and Horv\'ath (1997),
Gombay and Horv\'ath (1995, 2002)]  in that here we derive results assuming
 only $g(X_1)\in \hbox{DAN}$ and, consequently,
we may have $\sigma^2 = Eg^2(X_1) = \infty$.  The price we pay for this
 is the somewhat higher moment condition $E|h(X_1,X_2)|^{5/3} < \infty$ than that of the corresponding one with exponent 4/3 in Theorem 1. What is crucial in Theorem 2 in this regard is that the existence of the second moment of $h(X_1,X_2)$ is not assumed, for assuming the latter would exclude the possibility of having $\sigma^2=Eg^2(X_1)=\infty$ (cf. Remark 1).

\medskip
This paper is organized as follows. In the next section
we provide the proofs of main results.
Then, in Section \ref {sec3}, we investigate the asymptotic behaviour
of the $U-$type process $U_n(.)$ when it is based on kernels that are
antisymmetric, i.e.,  $h(x, y)$ in such that $h(x, y)=-h(y, x)$.
Throughout the paper $A, A_1,...$ will denote constants which may be different
in each appearance.

\section{Proofs of main results} \la {sec2}
We need some preliminaries to proving our main theorems.
The following lemma constitutes the key step. We note in passing that the three basic relations (11), (12), (13) of Lemma 1 are of interest on their own in studying U--statistics type processes, independently of their kernel function $h(\cdot, \cdot)$ being symmetric, or antisymmetric.

\begin{lem} \la {lem1} Let $\psi(x,y)$ be a measurable real valued symmetric
 function for which we have
\be
\int\psi(x, y) dF(y)=0 \la {rei0}
\ee and
 $E|\psi(X_1, X_2)|^{4/3}<\infty$. Then,
as $n\to\infty$,
 \be \f 1n \max_{1\le k\le
n-1}k^{-1/2}\Big|\sum_{i=1}^k\sum_{j=k+1}^n\psi(X_i, X_j)\Big| &=&
O_P(1),
\la {lem11}\\
\f 1{n} \max_{1\le k\le n-1}(n-k)^{-1/2}
\Big|\sum_{i=1}^k\sum_{j=k+1}^n\psi(X_i, X_j)\Big| &=&
O_P(1),\la {lem11a} \\
\f 1{n^{3/2}} \max_{1\le k\le n-1}\Big|\sum_{i=1}^k\sum_{j=k+1}^n\psi(X_i, X_j)\Big|
&=& o_P(1).
\la {lem12}
\ee
\end{lem}

{\it Proof.} We only prove (\ref{lem11}) and (\ref {lem12}). By virtue of
the symmetry of $\psi(x, y)$ and the i.i.d.\ properties of $X_i$,
the proof of (\ref {lem11a}) is similar to that of (\ref {lem11}). We omit the latter details.

 In order to  prove (\ref {lem11}), write
$$
\psi^*(x, y)=\psi(x, y)I_{\{|\psi(x, y)|\le
i^{3/2}\}}-\int\int\psi(u, v)I_{\{|\psi(u, v)|\le
i^{3/2}\}}dF(u)dF(v),
$$
$$
g^*(x)=\int \psi^*(x, y)dF(y) \hbox{ ~~and~~ } \psi^{**}(x, y)=\psi^*(x, y)-g^*(x)-g^*(y).
$$
It is readily seen that
\be
&&  E\big[\psi^*(X_i, X_j)\big] =0 \quad \mbox{and} \quad
 E\big[\psi^{**}(X_i, X_j)\mid X_i\big]=0, \quad \mbox{for all $i\not=j$.}  \la {rei1}
\ee
Having $E \big[\psi(X_i, X_j)\mid X_i\big]=0$ by (\ref {rei0}), we also have
\be
g^*(X_i)&=&E\big[\psi^{*}(X_i, X_j)\mid X_i\big] \no\\
&=&
E \big[\psi(X_i, X_j)I_{\{|\psi(X_i, X_j)|>i^{3/2}\}}\mid X_i\big]-
E\big[\psi(X_1, X_2)I_{\{|\psi|\ge
i^{3/2}\}}\big]. \la {yy6}
\ee

We now turn to the proof of (\ref {lem11}).  We have
\be
\f 1n \max_{1\le k\le n-1}k^{-1/2}\Big|\sum_{i=1}^k\sum_{j=k+1}^n\psi(X_i, X_j)\Big|
&\le & I_1(n)+I_2(n)+I_3(n), \la {ag1}
\ee
where
\bestar
I_1(n) &=&
\f 1n \max_{1\le k\le n-1}k^{-1/2}\Big|\sum_{i=1}^k\sum_{j\ne i\atop
j=1}^k\psi^{**}(X_i, X_j)\Big|, \no\\
 I_2(n) &=&
\f 1n \max_{1\le k\le n-1}k^{-1/2}\Big|\sum_{i=1}^k\sum_{j=1\atop
j\ne i}^n\psi^{**}(X_i, X_j)\Big|,\no\\
I_3(n) &=&
\f 1n \max_{1\le k\le n-1}k^{-1/2}\Big|\sum_{i=1}^k\sum_{j=k+1}^n
(\psi(X_i, X_j)-\psi^*(X_i, X_j)+g^*(X_i)+g^*(X_j))\Big|.
\eestar

We next prove $I_t(n)=O_P(1)$ for $t=1,2,3$ and then (\ref {lem11}) follows accordingly.

First consider $t=1$. Write $Y_i=\sum_{j=1}^{i-1}\psi^{**}(X_i, X_j)$.
Note that $E(Y_iY_k)=0$ for all $i\not= k$ by (\ref {rei1}).
It is readily  seen that
\be
E\Big|\sum_{i=2}^{\infty} i^{-3/2}\, Y_i\Big|^2 &= &
\sum_{i=2}^{\infty} i^{-3}\, EY_i^2 \le A\, \sum_{i=2}^{\infty} i^{-2} \,
E\big[\psi^2(X_1, X_2)I_{|\psi|\le i^{3/2}}\big]\no\\
&\le &A\,\sum_{k=1}^{\infty} E\big[\psi^2(X_1, X_2)
I_{(k-1)^{3/2}<|\psi|\le k^{3/2}}\big]\sum_{i=k}^{\infty} i^{-2} \no\\
&\le &A\, E|\psi(X_1, X_2)|^{4/3} <\infty. \la {yy1}
\ee
This, together with  the Kronecker lemma, implies that
$k^{-3/2}\sum_{i=1}^kY_i\to 0$, a.s.,
and hence $I_1(n)=O_P(1)$, since
 $I_1(n)\le 2\,\max_{1\le k\le n-1}k^{-3/2}
\big|\sum_{i=2}^k Y_i\big|$.

Secondly we  prove $I_2(n)=O_P(1)$. Write $Z_{in}=\sum_{j=1\atop
j\ne i}^n\psi^{**}(X_i, X_j) $.
By noting that, for any $a_i$ and $k\ge 1$,
$
\f 1k \sum_{i=1}^ka_i= b_k-\f 1k \sum_{i=1}^{k-1}b_i,
$
where $b_i=\sum_{t=1}^ia_t/t$, it follows that
\be
I_2(n) \le \f 1{n^{1/2}} \max_{1\le k\le n-1}k^{-1}\Big|\sum_{i=1}^kZ_{in}
\Big|\le \f 2{n^{1/2}}
\max_{1\le k\le n-1}\Big|\sum_{i=1}^k \f 1iZ_{in}\Big|. \la {yy2}
\ee
Therefore, it only needs to be  shown that, uniformly in $n\ge 1$,
\be
\f 1{n^{1/2}} E\Big|\sum_{i=1}^{\infty}\f 1{i}Z_{in}\Big|\le A<\infty.
\la {lem13}
\ee
Indeed, the result (\ref {lem13}) implies that
$\f 1{n^{1/2}} \Big|\sum_{i=1}^{\infty}\f 1{i}Z_{in}\Big|\le A<\infty$ {\it a.s.},  and
$\f 1{n^{1/2}} \Big|\sum_{i=k}^{\infty}\f 1{i}Z_{in}\Big|\to 0, a.s.,$ as $k\to\infty$,
 uniformly in $n\ge 1$.
This, together with (\ref {yy2}), yields
$$
I_{2n}\le \f 2{n^{1/2}}
\max_{1\le k\le N}\Big|\sum_{i=1}^{k} \f 1iZ_{in}\Big|\le
2A+ \max_{k\ge 1}\f 2{n^{1/2}} \Big|\sum_{i=k}^{\infty}\f 1{i}Z_{in}\Big|=O_P(1).
$$
The proof of (\ref {lem13}) follows from a similar argument
as in the proof of (\ref {yy1}). In fact,  for all $n\ge 1$, we have
\bestar
\f 1{n^{1/2}} E\Big|\sum_{i=1}^{\infty}\f 1{i }Z_{in}\Big|
&\le & \f 1{n^{1/2}} \Big[ E\Big|\sum_{i=1}^{\infty}\f 1{i }Z_{in}\Big|^2\Big]^{1/2} \no\\
&=&\f 1{n^{1/2}} \Big[\sum_{i=1}^{\infty}\f 1{i^2} E\Big(Z_{in}\Big)^2\Big]^{1/2}\\
&\le&A\,\Big[\sum_{i=1}^{\infty}\f 1{i^2}
E\psi^2(X_1,X_2)I_{(|\psi|\le i^{3/2})}\Big]^{1/2}\\
&< & A \big[E|\psi(X_1, X_2)|^{4/3}\big]^{1/2}< \infty,
\eestar
which yields (\ref {lem13}).

Finally we prove $I_3(n)=O_P(1)$.
Recalling (\ref {yy6}) and  $E\psi(X_1, X_2)=0$, we have
\bestar
\Lambda_{i,j}&:=& |\psi(X_i, X_j)-\psi^*(X_i, X_j)+g^*(X_i)+g^*(X_j)|\no\\
&\le&
|\psi(X_i, X_j)|I_{|\psi|\ge i^{3/2}}\ +\
E\big[|\psi(X_i,X_j)|I_{(|\psi|>i^{3/2})}|X_i\big]\no\\
&&\ +\
E\big[|\psi(X_i,X_j)|I_{(|\psi|>i^{3/2})}|X_j\big]\
+\ E\big[|\psi(X_1,X_2)|I_{(|\psi|>i^{3/2})}\big].
\eestar
This implies that $E(\Lambda_{i,j})\le 4 E\big[|\psi(X_1,X_2)|I_{(|\psi|>i^{3/2})}$, and hence
\be
E I_3(n)
&\le&
\f 1n E\Big[\max_{1\le k\le n-1}k^{-1/2}\sum_{i=1}^k\sum_{j=1\atop j\not =i}^n\Lambda_{i,j}\Big]\no\\
&\le & \f 1n \sum_{i=1}^{\infty}{i^{-1/2}}
\sum_{j=1\atop j\not =i}^nE(\Lambda_{i,j})\no\\
& \le&
4\sum_{i=1}^{\infty} \f 1{i^{1/2}} \,
E\big[|\psi(X_1,X_2)|I_{(|\psi|>i^{3/2})}\big] \no\\
& \le & \sum_{k=1}^{\infty} \,
E\big[|\psi(X_1,X_2)|I_{(k^{3/2}<|\psi|\le (k+1)^{3/2})}\big]\,
 \sum_{i=1}^{k} \f 1{i^{1/2}} \no\\
&\le &  A\, E |\psi(X_1,X_2)|^{4/3} <\infty, \la {yy8}
\ee
uniformly for all $n\ge 1$. By Markov's inequality, we obtain $I_3(n)=O_P(1)$.
The proof of (\ref {lem11}) is now complete.

The proof of (\ref {lem12}) is similar to that of (\ref {lem11}), but we have to use a
different truncation. In the following, we let
$$
\psi^*(x, y)=\psi(x, y)I_{\{|\psi(x, y)|\le
n^{3/2}\}}-\int\int\psi(u, v)I_{\{|\psi(u, v)|\le
n^{3/2}\}}dF(u)dF(v),
$$
$g^*(x)=\int \psi^*(x, y)dF(y)$ and $\psi^{**}(x, y)=\psi^*(x, y)-g^*(x)-g^*(y).$
It follows easily that
\be
\f 1{n^{3/2}} \max_{1\le k\le n-1}\, \Big|\sum_{i=1}^k\sum_{j=k+1}^n\psi(X_i, X_j)\Big|
&\le & \f 12\, \Big[I_0^*(n)+I_1^*(n)+I_2^*(n)\Big]+I_3^*(n), \la {ag1a}
\ee
where $I_0^*(n)=\f 1{n^{3/2}} \Big|\sum_{i=1}^n\sum_{j\ne i\atop
j=1}^n\psi^{**}(X_i, X_j)\Big|,$
\bestar
I_1^*(n) &=&
\f 1{n^{3/2}} \,\max_{1\le k\le n-1}\Big|\sum_{i=1}^k\sum_{j\ne i\atop
j=1}^k\psi^{**}(X_i, X_j)\Big|, \no\\
 I_2^*(n) &=&
\f 1{n^{3/2}}  \max_{1\le k\le n-1}\Big|\sum_{i=k+1}^n\sum_{j=k+1\atop
j\ne i}^n\psi^{**}(X_i, X_j)\Big|,\no\\
I_3^*(n) &=&
\f 1{n^{3/2}}  \max_{1\le k\le n-1}\Big|\sum_{i=1}^k\sum_{j=k+1}^n
(\psi(X_i, X_j)-\psi^*(X_i, X_j)+g^*(X_i)+g^*(X_j))\Big|.
\eestar
It is readily seen that
\bestar
E\big[I_0^*(n)\big]^2 &\le& A\, n^{-1}\,E\psi^2(X_1, X_2)I_{|\psi|\le n^{3/2}} \no\\
 &\le &A\, \Big[\ep^{-2}\, n^{-1/3}\, E|\psi(X_1, X_2)|^{4/3}
 +E|\psi(X_1, X_2)|^{4/3}I_{|\psi|\ge n}\Big] \no\\
 &\to& 0, \quad \mbox{as $n\to\infty$.}
\eestar
This yields $I_0^*(n)=o_P(1)$.
Noting that $\{\sum_{j=2}^kY_j, {\cal F}_k, 2\le k\le n\}$ is a martingale,
 where $Y_j=\sum_{i=1}^{j-1}\psi^{**}(X_i, X_j)$ and ${\cal F}_k=\si\{X_1, ...,X_k\}$,
 it follows from the well-known Maximum inequality for martingales that, for any $\ep>0$,
 \bestar
 P(I_1^*(n)\ge \ep) &\le& 4\ep^{-2}\,n^{-3}\,
 E\max_{1\le k\le n-1}\big|\sum_{j=2}^kY_j\big|^2
 \le A\,\ep^{-2}\,n^{-3}\,\sum_{j=2}^n EY_j^2 \no\\
 &\le& A\,\ep^{-2}\, n^{-1}\, E\psi^2(X_1, X_2)I_{|\psi|\le n^{3/2}} \no\\
 &\le &A\, \Big[\ep^{-2}\, n^{-1/3}\, E|\psi(X_1, X_2)|^{4/3}
 +E|\psi(X_1, X_2)|^{4/3}I_{|\psi|\ge n}\Big] \no\\
 &\to& 0, \quad \mbox{as $n\to\infty$.}
 \eestar
 This yields $I_1^*(n)=o_P(1)$. By a similar argument as in the proof for  $I_1^*(n)=o_P(1)$,
 we have $I_2^*(n)=o_P(1)$.
 As for $I_3^*(n)$, by using a similar argument as in the proof of (\ref {yy8}),
we obtain
 \bestar
 E \big|I_3^*(n)\big| &\le &\,\f 1{n^{3/2}} \,\sum_{i=1}^n\sum_{j=1\atop j\not= i}^n
E\,\Big|\psi(X_i, X_j)-\psi^*(X_i, X_j)+g^*(X_i)+g^*(X_j)\Big| \no\\
&\le& 4\, n^{1/2}\, E\big[|\psi(X_1, X_2)|I_{|\psi|\ge n^{3/2}}\big] \no\\
&\le& 4\, E\big[|\psi(X_1, X_2)|^{4/3}I_{|\psi|\ge n^{3/2}}\big] \to 0,
 \eestar
as $n\to\infty$, which implies that $I_3^*(n)=o_P(1)$. Taking all the respective estimates for
$I_{t}^{*}(n), t=0,1,2,3$ into (\ref {ag1a}), we obtain the required (\ref {lem12}).
The proof of Lemma \ref {lem1}
is now complete.

\medskip
The next two lemmas are due to CsCsHM (1986)
[cf. Lemma A.5.1
and Theorem A.5.1  respectively in Cs\"og\H{o} and Horv\'ath (1997)].
Proofs of Lemmas \ref {lem3} and \ref {alem3} can also be found in
Section 4.1 of Cs\"org\H{o} and Horv\'ath (1993).

\begin{lem} \la {lem3} Let $q(t)\in Q$.  If $I(q,c)<\infty$ for some $c>0$,
 then
\bestar
\lim_{t\downarrow 0}\,t^{1/2}/q(t)=0\quad\mbox{and}\quad
\lim_{t\uparrow 1}\,(1-t)^{1/2}/q(t)=0.
\eestar

\end{lem}

\begin{lem} \la {alem3} Let  $\{W(t), 0\le t<\infty\}$
be a standard Wiener process and  $q(t)\in Q$. Then,

{\rm (a)} $I(q,c)<\infty$ for all $c>0$ if and only if
\bestar
\limsup_{t\downarrow 0}\,|W(t)|/q(t)=0,\ a.s.\ \ \mbox{and}\ \
\limsup_{t\uparrow 1}\,|W(1)-W(t)|/q(t)=0, \ a.s.
\eestar

{\rm (b)} $I(q,c)<\infty$ for some $c>0$ if and only if
\bestar
\limsup_{t\downarrow 0}\,|W(t)|/q(t)<\infty,\
 a.s. \ \ \mbox{and}\ \
\limsup_{t\uparrow 1}\,|W(1)-W(t)|/q(t)<\infty,
\ a.s.
\eestar
\end{lem}

\bigskip
We are now ready to prove our main theorems.

\medskip
\textbf{Proof of Theorem \ref {th1}}. Together with the notation
 as in Section \ref {sec1}, we write $\psi(x, y)=h(x, y)-\theta-g(x)-g(y)$
and $T_n(t)=W_{[(n+1)t]}, 0\le t\le 1$, where
$$
W_k= (n-k)\sum_{j=1}^kg(X_j)+k\sum_{j=k+1}^ng(X_j).
$$
Noting that $g(X_j)$ are i.i.d.\ random variables with $Eg(X_1)=0$ and $\si^2=Eg^2(X_1)<\infty$,
along the lines of the proof of (2.1.45) in Cs\"org\H{o} and Horv\'ath  (1997),
 on an appropriate probability space for $X,X_1, X_2,\cdots$
we can define a sequence of Gaussian processes $\{\Gamma_n(t), 0\le t\le 1\}$
such that, for each $n\ge 1$,
$$
\{\Gamma_n(t), 0\le t\le 1\}{=}_d
\{\Gamma(t), 0\le t\le 1\},
$$
and if  $q\in Q$ and $I(q,c)<\infty$ for some $c>0$, then,
as $n\to\infty$,
\be
\sup_{1/n\le t\le (n-1)/n}\Big|\,n^{-3/2}\sigma^{-1}T_n(t)-
\Gamma_n(t)\,\Big|\Big/q(t)&=&o_P(1).
\la {th11ag}
\ee
By virtue of (\ref {th11ag}), Theorem \ref{th1} will follow if we prove
\be
J_n &:=&\sup_{1/n\le t\le (n-1)/n}\Big|n^{-3/2}U_n(t)-n^{-3/2}T_n(t)\Big|\Big/q(t)\ =\ o_P(1).
\la {try1}
\ee

In order to prove (\ref {try1}), write $V_n(t)=W_{[(n+1)t]}^*$, where
$
W_k^*= \sum_{j=1}^k\sum_{j=k+1}^n\psi(X_i, X_j).
$
Note that $E\big(\psi(X_1, X_2)\mid X_1\big)=E\big(\psi(X_1, X_2)\mid X_2\big)=0$ and
$$E|\psi(X_1, X_2)|^{4/3}\ \le\ A\, E|h(X_1, X_2)|^{4/3}\ <\ \infty.$$
 It follows  from (\ref {lem12})
that
\bestar
J_n^{(1)} &:=&\sup_{\delta\le t\le 1-\delta}
\big|n^{-3/2}V_n(t)\big|\Big/q(t) \no\\
&\le& \f 1{n^{3/2}} \max_{1\le k\le n-1}\Big|\sum_{i=1}^k\sum_{j=k+1}^n\psi(X_i, X_j)\Big| \,
\sup _{\delta\le t\le 1-\delta}\,q^{-1}(t)
= o_P(1),
\eestar
for any $\delta\in (0,1)$ and $q\in Q$. Let $\delta>0$ be so small that
$q(t)$ is already  nondecreasing   on $(0,\delta)$
and nonincreasing on $(1-\delta, 1)$ and let $n$ be so large such that $1/n\le \delta$.
It follows from (\ref {lem11}) and Lemma \ref {lem3} that
\bestar
J_n^{(2)} &:=&\sup_{0< t\le \delta}
\big|n^{-3/2}V_n(t)\big|\Big/q(t) \no\\
&\le& \f 1{n} \max_{1\le k\le n-1}k^{-1/2}\Big|\sum_{i=1}^k\sum_{j=k+1}^n\psi(X_i, X_j)\Big| \,
\sup _{0< t\le \delta}\,t^{1/2}/q(t)
= o_P(1),
\eestar
when $n\to\infty$ and then $\delta\to 0$. Similarly,
we have also
\bestar
J_n^{(3)} &:=&\sup_{1-\delta\le  t<1}
\big|n^{-3/2}V_n(t)\big|\Big/q(t) \no\\
&\le& \f 1{n} \max_{1\le k\le n-1}(n-k)^{-1/2}\Big|\sum_{i=1}^k\sum_{j=k+1}^n\psi(X_i, X_j)\Big| \,
\sup _{1-\delta\le  t<1}\,(1-t)^{1/2}/q(t)\no\\
&=&
 o_P(1),
\eestar
when $n\to\infty$ and then $\delta\to 0$. By virtue of these estimates, it is readily seen that
\be
J_n &\le &J_n^{(1)}+J_n^{(2)}+ J_n^{(3)}+
A\, n^{-1/2}\,\sup_{1/n\le t\le (n-1)/n}1/q(t)=o_P(1),
\ee
which yields (\ref {try1}). The proof of Theorem \ref {th1} is now complete.

\bigskip
\textbf{Proof of Corollary \ref {cor1}.} Having Theorem \ref {th1},
Lemmas \ref {lem3}-\ref{alem3} and the result (\ref {try1}),
  the proof of Corollary \ref {cor1} is the same as that given in the proof of Theorem 2.4.2
in  Cs\"org\H{o} and Horv\'ath (1997), and hence the details are omitted.

\bigskip

\textbf{Proof of Theorem \ref {th2}.}  We first prove (\ref {th15}).
 It is readily seen that
\be
\hat {U_n}(t)
&
=&n^{-3/2}(\hat {\sigma})^{-1}
\big\{Z_{[(n+1)t]}-n^2t(1-t) \theta\big\}
+t(1-t)n^{1/2}(\hat {\sigma})^{-1}(\hat {\theta}-\theta)\no\\
&=&
\left\{\f {\sum_{j=1}^ng^2(X_j)}{n\hat {\sigma}^2}\right\}^{1/2}\,
n^{-1} \Big\{\sum_{j=1}^ng^2(X_j)
\Big\}^{-1/2}U_n(t)+t(1-t)n^{1/2}(\hat {\sigma})^{-1}(\hat {\theta}-\theta).\no\\
\la {ah1}
\ee
Furthermore  $U_n(t)=T_n(t)+V_n(t)$,
where $T_n(t)$ and $V_n(t)$ are defined as in the proof of Theorem \ref {th1}.
Recalling that $g(X_1)$ is in the domain of attraction of the normal law,
as in the proof of Theorem 5.2 of Cs\"org\H{o}, Szyszkowicz and  Wang [CsSzW] (2004) with minor modifications, we have that
on an appropriate probability space for $X,X_1, X_2,\cdots,$
we can define a sequence of Gaussian processes $\{\Gamma_n(t), 0\le t\le 1\}$
such that (\ref {ad10}) holds true, and as $n\to\infty$,
\bestar
\sup_{0< t<1}\Big|\,n^{-1} \Big\{\sum_{j=1}^ng^2(X_j)
\Big\}^{-1/2}\,T_n(t)-
\Gamma_n(t)\,\Big|\Big/q(t)&=&o_P(1),
\eestar
if and only if  $I(q,c)<\infty$ for all $c>0$. Therefore, to prove (\ref {th15}),
 it suffices to show that
\be
n^{-1} \Big\{\sum_{j=1}^ng^2(X_j)
\Big\}^{-1/2} \, \sup_{0<t<1}|V_n(t)|/q(t) &=&o_P(1), \la {th78} \\
\Big\{n^{-1} \sum_{j=1}^ng^2(X_j)\Big\}^{-1}\,\hat {\sigma}^2-1
&=&o_P(1),
\la {th22}
\ee
and
\be
\qquad \qquad n^{1/2}(\hat {\sigma})^{-1}(\hat {\theta}-\theta)
&=&o_P(1).
\la {th23}
\ee
The proof of (\ref {th78}) is simple and in fact (\ref {th78}) holds true
if $q(x)$ satisfies $I(q,c)<\infty$ for some $c>0$. Indeed, since $g(X_1)$
is in the domain of attraction of the normal law, we have
$\f 1{b_n}\sum_{j=1}^n g^2(X_j)\to_P 1$, where $b_n=n\,l(n)$
with that $l(n)=Eg^2(X_1)$ if $Eg^2(X_1)<\infty$ or $l(n)\to \infty$ if $Eg^2(X_1)=\infty$.
On the other hand, as in the proof of (\ref {try1}),
 $n^{-3/2}\sup_{0<t<1}|V_n(t)|/q(t)=o_P(1)$ even when $q(x)$ satisfies $I(q,c)<\infty$ for some $c>0$, and hence
 (\ref {th78}) follows immediately from these facts.

We next prove (\ref {th22}). The claim (\ref {th23}) follows  by using (\ref {th22}),
and hence the details are omitted.
Without loss of generality, we  assume $\theta=0$. We may rewrite $\hat {\sigma}^2$\
as
\bestar
\hat{\sigma}^2&=&\frac 1{n(n-1)^2}\sum_{i\ne j\ne k}h(X_i, X_j)h(X_i, X_k)
+\frac 1{n(n-1)^2}\sum_{i\ne j}h^2(X_i, X_j)-\hat {\theta}^2\\
&:=& W_{n1}+W_{n2}-\hat {\theta}^2.
\eestar
Recalling $E|h(X_1,X_2)|^{5/3}<\infty$, it follows  from
a Marcinkiewicz type strong  law for $U$-statistics that
$W_{n2}-\hat {\theta}^2\to 0, a.s.$ [see Gine and Zinn (1992), for example]. Therefore
(\ref {th22}) will follow if we  prove
 \be
\Big\{n^{-1} \sum_{j=1}^ng^2(X_j)\Big\}^{-1}\,W_{n1}-1
&=&o_P(1).
\la {th22a}
\ee
Write, for $i\ne j\ne k$,
\bestar
h_{ij}^{(1)}
&=&h(X_i, X_j)I_{(|h|\le n^{6/5})}, \qquad g^{(1)}(X_i)\ = \
E\big(h_{ij}^{(1)}\big|X_i\big),\\
\psi_{ijk} &=&   h_{ij}^{(1)}\, h_{ik}^{(1)}-Eh_{ij}^{(1)}\, h_{ik}^{(1)},\\
\varphi_i^{(1)} &=& E\big(\psi_{ijk}\big|X_i\big), \quad
\varphi_j^{(2)} \ =\ E\big(\psi_{ijk}\big|X_j\big), \quad
\varphi_k^{(3)} \ =\ E\big(\psi_{ijk}\big|X_k\big).
\eestar
Noting that
$E\big\{h_{ij}^{(1)}\, h_{ik}^{(1)}\big|X_i\big\}=\big\{g^{(1)}(X_i)\big\}^2$,
it is readily seen that $\varphi_i^{(1)}=\big\{g^{(1)}(X_i)\big\}^2-
E\big[h_{ij}^{(1)}\, h_{ik}^{(1)}\big]$, and
\bestar
\sum_{i\ne j\ne k}h_{ij}^{(1)}\, h_{ik}^{(1)}
&=&\sum_{i\ne j\ne k}\psi_{ijk}+
\sum_{i\ne j\ne k}E\big[h_{ij}^{(1)}\, h_{ik}^{(1)}\big] \\
&=&\sum_{i\ne j\ne k}\big\{g^{(1)}(X_i)\big\}^2+\sum_{i\ne j\ne k}
\big\{\varphi_j^{(2)}+\varphi_k^{(3)}\big\}\\
&&\quad +\ \sum_{i\ne j\ne k}\big(\psi_{ijk}-\varphi_i^{(1)}
-\varphi_j^{(2)}-\varphi_k^{(3)}\big)\\
&:=& V_{n1}+V_{n2}+V_{n3}.
\eestar
In the next paragraph, we will show that
\be
\Big\{n^{-1} \sum_{j=1}^ng^2(X_j)\Big\}^{-1}\,\Big(n^{-3}\,V_{n1}\Big)-1
&=&o_P(1),
\la {who}\\
n^{-3}\, \big(V_{n2}+V_{n3}\big)&=&o_P(1).
\la {who1}
\ee
It follows from (\ref {who}) and (\ref {who1}) that
\be
\Big\{n^{-1} \sum_{j=1}^ng^2(X_j)\Big\}^{-1}\,n^{-3}
\sum_{i\ne j\ne k}h_{ij}^{(1)}\, h_{ik}^{(1)}\ -\ 1
&=&o_P(1),
\la {who2}
\ee
and then (\ref {th22a}) follows  from (\ref {who2}) and
\bestar
P\(\sum_{i\ne j\ne k}h_{ij}\, h_{ik}
\neq\sum_{i\ne j\ne k}h_{ij}^{(1)}\, h_{ik}^{(1)} \) &\leq& n^2\,
P\big(|h(X_1, X_2)|\geq n^{6/5}\big) \no\\
&\le& E|h(X_1, X_2)|^{5/3}I_{|h|\ge n^{6/5}}\to 0.
\eestar

We are to prove (\ref {who}) and (\ref {who1}) now.  Consider (\ref {who}) first.
By noting that $g^{(1)}(X_1)=g(X_1)-g^*(X_j)$,
where $g^*(X_j)=E\big\{h(X_1,X_2)I_{(|h|\ge n^{6/5})}\big|X_1\big\}$, we have
\bestar
\Big|\sum_{j=1}^n\Big[\big\{g^{(1)}(X_j)\big\}^2-g^2(X_j)\Big]\Big| &\le&
\sum_{j=1}^n\Big[2|g(X_j)|\, |g^*(X_j)|+|g^*(X_j)|^2\Big] \no\\
&\le& 2 \Big[\sum_{j=1}^ng^2(X_j)\Big]^{1/2}\,\Big[\sum_{j=1}^n\{g^*(X_j)\}^2\Big]^{1/2}+
\sum_{j=1}^n\{g^*(X_j)\}^2.
\eestar
Now, since $g(X_1)$ is in the domain of attraction of the normal law
[which implies that $\f 1{n}\sum_{j=1}^ng^2(X_j)\to_P C>0$, where $C$ may be $\infty$],
simple calculations show that
(\ref {who}) will follow if we prove
\be
\frac 1n\sum_{j=1}^n\{g^*(X_j)\}^2=o_P(1).\la {ag20}
\ee
In fact, for any $\ep>0$, we have
\bestar
P\Big(\sum_{j=1}^n\{g^*(X_j)\}^2\ge \ep\, n\Big) &\le & \ep^{-1/2}n^{-1/2}\sum_{j=1}^n
E|g^*(X_j)| \no\\
&\le& \ep^{-1/2}n^{1/2} E|h(X_1,X_2)|I_{(|h|\ge n^{6/5})}\no\\
&\le& \ep^{-1/2}\,
E|h(X_1,X_2)|^{5/3}I_{(|h|\ge n^{6/5})}
\to 0,
\eestar
as $n\to\infty$. This implies (\ref {ag20}) and hence completes the proof of (\ref {who}).

We next prove (\ref {who1}). By noting that $n^{-3}V_{n3}$
is a degenerate $U$-statistic of order $3$, it follows from  moment inequality
 for degenerate $U$-statistics
(see, Borovskikh (1996), for example)  that, for any $\ep>0$,
\be
P\big(|V_{n3}|\ge \ep n^3\big) &\le& \ep^{-5/3}\, n^{-5}\,
E|V_{n3}|^{5/3} \no\\
&\le&
A\,\ep^{-5/3}\, n^{-2}\, E\Big|\psi_{123}-\varphi_1^{(1)}-\varphi_2^{(2)}-\varphi_3^{(3)}\Big|^{5/3}
\no\\
&\le& A\, \ep^{-5/3}\, n^{-2}\, E|h(X_1,X_2)|^{10/3}I_{(|h|\leq n^{6/5})} \no\\
&\le & A\, \ep^{-5/3}\,\Big[n^{-1/3}+ E|h(X_1, X_2)|^{5/3}I_{(|h|\geq n^{1/2})}\Big]\to 0,
\la {who4}
\ee
as $n\to \infty$.
On the other hand, by noting that
\bestar
E\Big\{E\Big[h_{12}^{(1)}\, h_{13}^{(1)}\big|X_2\Big]\Big\}^2
&=&E\Big\{h_{12}^{(1)}\, h_{13}^{(1)}\,
E\Big[h_{42}^{(1)}\, h_{45}^{(1)}\big|X_2\Big]\Big\}\\
&=&E\Big[h_{12}^{(1)}\, h_{13}^{(1)}\,
h_{42}^{(1)}\, h_{45}^{(1)}\Big]\\
&\le & \Big[Eh^2(X_1, X_2)I_{|h|\le n^{6/5}}\Big]^2\le n^{4/5}\,
\Big\{E|h(X_1,X_2)|^{5/3}\Big\}^2,
\eestar
it is readily seen that, for any $\ep>0$,
\be
P\big(|V_{n3}|\ge \ep n^3\big) &\le& \ep^{-2}\, E\Big(n^{-3}V_{n2}\Big)^2 \no\\
&\le & A \,\ep^{-2}\,n^{-1} E\Big(\varphi_1^{(2)}+\varphi_1^{(3)}\Big)^2\no\\
&\le & A\, \ep^{-2}\,n^{-1} \[
 E\Big\{E\Big(h_{12}^{(1)}\, h_{13}^{(1)}\big|X_2\Big)\Big\}^2+
 \Big(E\Big\{h_{12}^{(1)}\Big\}^2\Big)^2\]\no\\
 &\le & A\,\ep^{-2}\, n^{-1/5}\,\Big\{E|h(X_1,X_2)|^{5/3}\Big\}^2 \to 0, \la {who5}
\ee
as $n\to \infty$. By virtue of  (\ref {who4}) and (\ref {who5}), we obtain (\ref {who1}).
The proof of (\ref {th15}) is now complete.

The result (\ref {adc12}) is a direct consequence of (\ref {th15}).
As for (\ref {th17}), by virtue of (\ref {ah1})-(\ref {th23})
(recalling that (\ref {th78}) still holds true for $q(x)$ satisfying
$I(q, c)<\infty$ for some $c>0$, as explained in its proof), it suffices to show that
\be
\sup_{0< t<1}\Big|\,n^{-1} \Big\{\sum_{j=1}^ng^2(X_j)
\Big\}^{-1/2}\,T_n(t)\Big|&\rightarrow_d& \sup_{0<t<1}|\Gamma(t)|/q(t)
\la {th17a}
\ee
if and only if $I(q,c)<\infty$ for some $c>0$, where
$T_n(t)=W_{[(n+1)t]}, 0\le t\le 1$, with
$$
W_k= (n-k)\sum_{j=1}^kg(X_j)+k\sum_{j=k+1}^ng(X_j).
$$
This follows from the same arguments as in the proof of Corollary 5.2 in CsSzW (2004),
and hence the details are omitted. This also completes the proof of Theorem \ref {th2}.

\section{Antisymmetric kernel} \la {sec3}

In this section we consider the asymptotics of  $U$-type processes with  antisymmetric kernel $h(x, y)$,
i.e., $h(x, y)=-h(y,x)$. This kind of kernels can not be symmetrized, but they
are especially useful to check the equality of distributions
for different groups of random variables since $\theta =Eh(X_1, X_2)=0$
whenever $X_1 =_d X_2$, if $E|h(X_1,X_2)|<\infty$.  Consequently, for antisymmetric kernels, $U_n(t) = Z_{[(n+1)t]}$ under $H_0$.. An example is given in Pettitt (1979),who used functions of the Mann-Whitney type statistics
\bestar
(12)^{1/2}n^{-3/2}\sum_{1\le i\le nt}\sum_{nt<j\le n}\, \mbox{sign}\, (X_i-X_j)
\eestar
to detect possible changes in distribution. Another important example is given by taking $H(x,y)=x-y$ for studying the probable error of a change in a mean.  We will say more about that in Remark 5.

 For the anti-symmetric kernel $h(x, y)$, by letting $g(t)=Eh(X_1, t)$, i.e., keeping our earlier notation with $\theta = 0$, 
we may write
\bestar
Z_k=\sum_{i=1}^k\sum_{j=k+1}^n\psi(X_i,X_j)+ n\,\Big[\sum_{i=1}^kg(X_i)-\f kn\,
\sum_{i=1}^ng(X_i)\Big],
\eestar
where $\psi(x, y)=h(x, y)+g(x)-g(y)$ with
$$
E\[\psi(X_1, X_2)\mid X_1\]=E\[\psi(X_1, X_2)\mid X_2\]=0.
$$

Since Lemma \ref{lem1} does not depend on the symmetry of the kernel,
similarly to the proofs of Theorems \ref {th1} and \ref {th2},
we have the following results for
$U$-type processes with  antisymmetric kernel $h(x, y)$,
which improve and generalize the similar earlier results of
Cs\"org\H{o} and Horv\'ath (1988a,\,b),
Szyszkowicz (1991, 1992) and those
 given in Section 2.4 of
 Cs\"org\H{o} and Horv\'ath (1997) along these lines.
 It is interesting  to note that
 the  Gaussian limit process that is shared by Theorems \ref {th1} and \ref {th2}
 and that shared by Theorems \ref {th3} and \ref {th4}
 are different,
  although they are of equal variance.
For further related  results, we refer to Janson and Wichura (1983), and
 Gombay (2000a, b, 2001, 2004).

 We continue to use the notations introduced in Section \ref {sec1}, but
  $U_n(t)$ and $\hat U_n(t)$ are now defined  in terms of antisymmetric kernel $h(x, y)=-h(y,x)$.
Consequently, under $H_0$, $\theta$ and $\hat\theta$ are both zero now.

\begin{th} \la {th3} Let $q\in Q$. Assume $H_0$, $0<\sigma^2<\infty$
and $E|h(X_1, X_2)|^{4/3}<\infty$.
Then, on an appropriate probability space for $X,X_1, X_2,\cdots,$
we can define a sequence of Brownian bridges  $\{B_n(t), 0\le t\le 1\}$
such that if $I(q,c)<\infty$ for some $c>0$, then as $n\to \infty$,
\be
\sup_{1/n\le t\le (n-1)/n}\Big|\,n^{-3/2}\sigma^{-1}U_n(t)-
B_n(t)\,\Big|\Big/q(t)&=&o_P(1).
\la {th11a10}
\ee
Consequently,

\quad {\rm(a)} as $n\to\infty$,
\begin{equation}
\sup_{0<t<1} \big|n^{-3/2} \sigma^{-1}U_n(t) - B_n(t)|/q(t) = o_P(1)
\label{eq38}
\end{equation}
if and only if $I(q,c)<\infty$ for all $c>0$;

\quad {\rm(b)}  as $n\to\infty$,
\begin{equation}
n^{-3/2}\sigma^{-1}U_n(\cdot)\,\Big/q(\cdot) \Rightarrow
B(\cdot)\,\Big/q(\cdot) \hbox{ on } (D[0,1], {\mathcal D}, ||/q||)
\la {adcor12s}
\end{equation}
if and only if
$I(q,c)<\infty$ for all $c>0$;

\quad {\rm(c)} as $n\to\infty$,
\be
n^{-3/2}\sigma^{-1}\sup_{0<t<1}|U_n(t)|\Big/q(t)&\rightarrow_d&
\sup_{0<t<1}|B(t)|/q(t)
\la {th13s}
\ee
if and only if $I(q,c)<\infty$ for some $c>0$, where, in (b) and (c),
$\{B(t), 0\le t\le 1\}$ is a Brownian bridge.
\end{th}

Theorem \ref {th3} is to be compared to Szyszkowicz (1991, Theorem 2.1)
[cf. Theorem 2.4.1 in Cs\"org\H{o} and Horv\'ath (1997)].

\begin{th} \la {th4} Let $q\in Q$. Assume $H_0$, $E|h(X_1, X_2)|^{5/3}<\infty$
and that $g(X_1) \in \hbox{DAN}$.
Then, on an appropriate probability space for $X,X_1, X_2,\cdots,$
we can define a sequence of Brownian bridges  $\{B_n(t), 0\le t\le 1\}$
such that, as $n\to\infty$,
\be
\sup_{0< t<1}\Big|\,\hat {U_{n}}(t)-
B_n(t)\,\Big|\Big/q(t)&=&o_P(1),
\la {th15s}
\ee
if and only if  $I(q,c)<\infty$ for all $c>0$. Consequently, as $n\to\infty$,
\be
\hat {U_{n}}(\cdot)/q(\cdot)&\Rightarrow&
B(\cdot)\,\Big/q(\cdot),\quad on\ (D[0,1], {\mathcal D},||/q||)
\la {adc12s}
\ee
 if and only if
$I(q,c)<\infty$ for all $c>0$, where
$\{B(t), 0\le t\le 1\}$ is a Brownian bridge. Furthermore, as $n\to\infty$,  we also have
\be
\sup_{0<t<1}|\hat {U_n}(t)|\Big/q(t)&\rightarrow_d& \sup_{0<t<1}|B(t)|/q(t)
\la {th17s}
\ee
if and only if $I(q,c)<\infty$ for some $c>0$.
\end{th}

\medskip\noindent{\bf Remark 4}~ As compared to Theorem 3, where it is assumed that $0 < \sigma^2 = Eg^2(X_1)<\infty$, in Theorem 4 we assume only that $g(X_1)$ is in the domain of attraction of the normal law and, consequently, we may have $\sigma^2 = Eg^2(X_1) = \infty$, just like in Theorem 2 (cf. Remark 3).

\medskip\noindent{\bf Remark 5}~ On taking $h(x, y)=x-y$, Theorem \ref {th4} essentially extends
Corollary 2.1.1 of Cs\"org\H{o} and Horv\'ath (1997)
[cf. Theorem 5.1 in CsSzW (2004)] and rhymes with Theorem 5.2 and Corollaries 5.1 and 5.2 of CsSzW (2004)
[cf. also Theorem 2.1 and Corollaries 2.1 and 2.2 of CsSzW (2006)], where we study the problem of change in the mean
in DAN directly via Theorem 2 and Corollaries 3 and 4 of CsSzW (2007), quoting these results without proof for the sake of studying the probable error of a change in a mean in the domain of attraction of the normal law. In this regard our present Theorems 2 and 4 can be viewed in part as extensions of the initial scope of our research in CsSzW (2007) on weighted approximations of self-normalized partial sum processes to those of Studentized U-statistics type processes with symmetric and antisymmetric kernel functions $h(\cdot,\cdot)$, whose respective projections $g(X_1)$ are in DAN.



\vspace{1cm}
\begin{center}
{\large{\bf REFERENCES}}
\end{center}

\newcommand{\itemr}[1] {\par\indent\hbox to -29.0pt{#1\hfil} \hangindent=20.0pt
\hangafter=1\hskip 10pt}

\itemr{}
Borovskikh, Yu. V.(1996). {\it $U$-statistics in Banach spaces.}
VSP, Utrecht.

\itemr{}
Borovskikh, Yu. V. (2002). On the normal approximation of
$U$-statistics. {\it Theory Probab.\ Appl.}  {\bf 45},  406--423.

\itemr{} Cs\"org\H{o}, M.\ (2002).  A glimpse of the impact of P\'al Erd\H{o}s on probability and statistics.  \textit{The Canadian Journal of Statistics} {\bf 30}, 493--556.

\itemr{} Cs\"org\H{o}, M., Cs\"org\H{o}, S., Horv\'ath, L.\ and Mason, D.~(1986).
Weighted empirical and quantile processes, {\it Ann. Probab.} {\bf
14}, 31-85.

\itemr{} Cs\"org\H{o}, M.\ and  Horv\'ath, L.~(1988a).
 Invariance principles for changepoint problems.
 {\it J. Multivariate Anal.}  {\bf 27}, 151--168.

\itemr{} Cs\"org\H{o}, M.\ and  Horv\'ath, L.~(1988b). Nonparametric methods
for changepoint problems, In {\it Handbook of Statistics}, Elsevier Science
Publisher B.V.,  403-425, North-Holland, Amsterdam.

\itemr{} Cs\"org\H{o}, M.\ and  Horv\'ath, L.~(1993).
{\it Weighted Approximations in Probability and Statistics}, Wiley
Series in Probability and Mathematical Statistics: Probability and
Mathematical Statistics, Wiley, Chichester.

\itemr{} Cs\"org\H{o}, M., and  Horv\'ath, L.~(1997). {\it Limit Theorems
in Change-Point Analysis}, Wiley Series in Probability and
Mathematical Statistics: Probability and Mathematical Statistics,
Wiley, Chichester.

\itemr{} Cs\"org\H{o}, M., Norvai\v{s}a, R.\ and Szyszkowicz, B.~(1999).
Convergence of weighted partial sums when the limiting distribution
is not necessarily Radon, {\it Stochastic Process. Appl.} {\bf 81},
81-101.
\itemr{}
Cs\"org\H{o}, M., Szyszkowicz, B.\ and  Wang, Q.\ (2004). On weighted
approximations and strong limit theorems for self-normalized partial
sums processes. In {\it Asymptotic methods in stochastics,}  489--521,
Fields Inst. Commun. {\bf 44}, Amer.\ Math.\ Soc., Providence, RI.

\itemr{}
Cs\"org\H{o}, M., Szyszkowicz, B.\ and  Wang, Q.\ (2006). Change in the mean in the domain of attraction of the normal law.
{\it Austrian Journal of Statistics}, {\bf 35}, 93-103.

\itemr{}
Cs\"org\H{o}, M., Szyszkowicz, B.\ and  Wang, Q.\ (2007). Weighted approximations in $D[0,1]$ with applications to self-normalized partial sum processes.  {\it Preprint}.

\itemr{} Gine, E., and Zinn, J. (1992). Marcinkiewicz type laws of large numbers and
convergence of moments for U-statistics. In {\it Probability in
Banach Spaces} (R. Dudley, M.~Hahn and J.~Kuelbs, eds) {\bf 8}
273-291, Birkhauser, Boston.

\itemr{}
Gombay, E. (2000a).  Comparison of $U$-statistics in the
change-point problem and in sequential change detection. Endre Cs\'aki
65. {\it Period.\ Math.\ Hungar.} {\bf  41}, 157--166.

\itemr{}
Gombay, E. (2000b).  $U$-statistics for sequential change detection.
{\it Metrika} {\bf 52}, 133--145.

\itemr{}
Gombay, E. (2001). $U$-statistics for change under alternatives.
{\it J.\ Multivariate Anal.}  {\bf 78}, 139--158.

\itemr{}
Gombay, E. (2004). $U$-statistics in sequential tests and change
detection. Abraham Wald centennial celebration: invited papers. Part
II. {\it Sequential Anal.}  {\bf 23}, 257--274.

\itemr{}
Gombay, E. and  Horv\'ath, L. (1995). An application of $U$-statistics
to change-point analysis.  {\it Acta Sci.\ Math.} (Szeged) {\bf  60},
345--357.

\itemr{}
Gombay, E. and  Horv\'ath, L. (2002). Rates of convergence for
$U$-statistic processes and their bootstrapped versions. {\it  J.\
Statist.\ Plann.\ Inference} {\bf 102}, 247--272.

\itemr{}
Halmos, P.R.\ (1946).  The theory of unbiased estimation.  {\it Ann.\ Math.\ Statist.} {\bf 17}, 34--43.

\itemr{}
Hoeffding, W.\ (1948).  A class of statistics with asymptotically normal distribution.  {\it Ann.\ Math.\ Statist.} {\bf 19}, 293--325.
\itemr{}
Janson, S. and  Wichura, M. J. (1984).
 Invariance principles for stochastic area and related stochastic integrals.
{\it Stochastic Process.\ Appl.} {\bf   16},  71--84.
\itemr{}
Pettitt, A. N. (1979). A nonparametric approach to the change-point problem.
{\it J. Roy.\ Statist.\ Soc.\ Ser.\ C} {\bf   28}, 126--135.

\itemr{}
Shorack, G.R.\ and Wellner, J.A.\ (1986).  \textit{Empirical Processes with Applications to Statistics}, Wiley, New York.

\itemr{} Szyszkowicz, B.~(1991). Weighted stochastic processes under contiguous
alternatives,  {\it C.R. Math.\ Rep.\ Acad.\ Sci.\ Canada } {\bf 13},
211-216.

\itemr{} Szyszkowicz, B.~(1992). {\it Weak Convergence of Stochastic Processes in Weighted
Metrics and their Applications to Contiguous Changepoint Analysis}. Ph.\ D.
 Dissertation, Carleton University.

\itemr{} Szyszkowicz, B.~(1996). Weighted approximations of partial sum processes
in $D[0,\infty)$. I, {\it Studia Sci.\ Math.\ Hungar.} {\bf 31},
323-353.

\itemr{} Szyszkowicz, B.~(1997). Weighted approximations of partial sum processes
in $D[0,\infty)$. II, {\it Studia Sci.\ Math.\ Hungar.} {\bf 33},
305-320.

\vskip1cm

\noindent
Mikl\'os Cs\"org\H{o} \\
School of Mathematics and Statistics\\
Carleton University\\
1125 Colonel By Drive\\
Ottawa, ON ~Canada K1S 5B6\\
{\tt mcsorgo@math.carleton.ca}
\\[0.5cm]
Barbara Szyszkowicz \\
School of Mathematics and Statistics\\
Carleton University\\
1125 Colonel By Drive\\
Ottawa, ON ~Canada K1S 5B6\\
{\tt bszyszko@math.carleton.ca}
\\[0.5cm]
Qiying Wang \\
School of Mathematics and Statistics\\
University of Sydney\\
NSW 2006, Australia\\
{\tt qiying@maths.usyd.edu.au}

\end{document}